\documentclass[12pt]{amsart}

\usepackage{amssymb}
\usepackage{mathabx}
\usepackage{mathrsfs}

\usepackage{enumitem}

\makeatletter
\@namedef{subjclassname@2010}{%
  \textup{2010} Mathematics Subject Classification}
\makeatother

\usepackage[T1]{fontenc}

\newtheorem{theorem}{Theorem}[section]

\newtheorem{proposition}[theorem]{Proposition}

\theoremstyle{definition}
\newtheorem{definition}[theorem]{Definition}

\newtheorem{example}[theorem]{Example}

\numberwithin{equation}{section}

\begin{document}


\baselineskip=17pt


\title[Invariant submanifolds for affine control systems]{Invariant submanifolds for affine control systems}

\author[C.-K. Han]{Chong-Kyu Han}
\address[C.-K. Han]{Department of Mathematical Sciences\\ Seoul National University\\ 1 Gwanak-ro, Gwanak-gu\\ Seoul 08826, Republic of Korea}
\email{ckhan@snu.ac.kr}
%

\author[H. Kim]{Hyeseon Kim}
\address[H. Kim]{Research Institute of Mathematics\\ Seoul National University\\ 1 Gwanak-ro, Gwanak-gu\\ Seoul 08826, Republic of Korea}
\email{hop222@gmail.com}

\date{}

\begin{abstract}
Given an affine control system $\dot{\mathbf x} = f({\mathbf x}) + \sum_{j=1}^m  g_j({\mathbf x}) u_j$ we present a method  of construction of  submanifolds that are invariant under controls  assuming  that the linear span of $f,  g_1, \ldots, g_m$ has constant rank. We use the  method of reduction of Pfaffian systems to a largest integrable subsystem and  finding the first integrals and the generalized first integrals for the vector fields $f$ and $g_j$'s.

\end{abstract}



\maketitle

\section{Introduction and the statement of the main results}\label{intro}
Let $M$ be a connected smooth~($C^{\infty}$) manifold of dimension $n$ and $\mathcal{U}$ be a set of admissible controls $u: [0,\infty) \rightarrow\mathbb{R}^m$, $m\leq n$. We consider an \emph{affine control system,} which is a system of differential equations of the form:
\begin{equation}\label{affine_control_system}
\dot{x}:=\dfrac{dx}{dt}=f(x)+\sum_{j=1}^{m}g_{j}(x)u_{j},
\end{equation}where $x=(x_1,\ldots,x_n)$ are local coordinates of $M$, $f$ and $g_1,\ldots,g_m$ are smooth vector fields on $M$, and $u=(u_1,\ldots,u_m)\in\mathcal{U}$. Here the control $u(t)$ can be chosen variously, for instance, to be piecewise continuous, measurable, smooth, and so forth.   In this paper we shall assume  that $\mathcal U$ is the set of all piecewise constant functions  with finitely many discontinuities.

Our viewpoint is local,  thus $M$ can be regarded throughout this paper as an open ball of $\mathbb{R}^n$ centered at a reference point. 
For a point $p\in M$ let $x(t)$ be the solution of

\begin{equation*}\label{IVP}
\left\{
\begin{aligned}
&\dot{x}=f(x)+\sum_{j=1}^{m}g_{j}(x)u_{j}(t),\\
&x(0)=p.
\end{aligned}
\right.
\end{equation*}

The solution $x(t)$  is continuous,  piecewise smooth and uniquely determined by the choice of $p$ and $u(t)$, which we shall denote by $\Gamma(p,t,u)$.
\begin{definition}
A submanifold $N$ of $M$ is said to be \emph{invariant under controls} of \eqref{affine_control_system} if $x_{0}\in N$ implies that $\Gamma(x_0,t,u)\in N$ for all possible choices of $u\in \mathcal U $ and for all $t\ge 0.$
\end{definition}

The purpose of this paper is to present  a  systematic description and a  method of construction of invariant submanifolds for \eqref{affine_control_system}.  By a \emph{distribution} $\mathcal{D}$ of rank $k$ we shall mean a smooth sub-bundle $\mathcal{D}$ of fibre dimension $k$  of the tangent bundle $TM$. A smooth real-valued function $\rho$ is called a \emph{first integral} of $\mathcal{D}$ if
\begin{equation}\label{FI}
X\rho=0
\end{equation}for any section $X$ of $\mathcal{D}$. A distribution $\mathcal{D}$ is said to be \emph{integrable} if 
\[
[\mathcal{D},\mathcal{D}]\subset\mathcal{D},
\]
which means that for any sections $X$ and $Y$ of $\mathcal{D}$ their Lie bracket $[X,Y]$ is a section of $\mathcal{D}.$ ~~   A system $\rho=(\rho^1,\ldots,\rho^d)$ of smooth real-valued functions on $M$ is said to be \emph{non-degenerate} if
\[
d\rho^1\wedge\cdots\wedge d\rho^d\neq0  \quad\text{on a neighborhood of }  \{\rho=0\}.
\]

\begin{theorem}[Frobenius]\label{Frobenius}
Let $M$ be a smooth manifold of dimension $n$ and $\mathcal{D}$  a distribution of rank $k$. If $\mathcal{D}$ is integrable, then there exists a non-degenerate system of $n-k$ first integrals $\rho=(\rho^1,\ldots,\rho^{n-k})$. 
\end{theorem}Given a distribution $\mathcal{D}$ of rank $k \geq2, $  \eqref{FI} is over-determined, and there are no solutions generically. Theorem~\ref{Frobenius} is an extreme case that \eqref{FI} has $(n-k)$ independent solutions.  This is the largest possible dimension for the solution space.  If $\mathcal{D}$ is not integrable one might construct an integrable  distribution of the smallest rank that contains $\mathcal{D}$ as a sub-bundle as follows: Set 
\[\mathcal{D}^{(0)}=\mathcal{D}\]
and for each $j=1,2,\ldots$, 
\[
\mathcal{D}^{(j)}=\mathcal{D}^{(j-1)}+[\mathcal{D}^{(j-1)},\mathcal{D}^{(j-1)}],
\]assuming $\mathcal{D}^{(j-1)}$ and $\mathcal{D}^{(j)}$ have constant rank.
 Let $\nu$ be the smallest non-negative integer such that 
\[
\mathcal{D}^{(\nu)}=\mathcal{D}^{(\nu+1)}.
\]We call the  sequence
\begin{equation}\label{derived_flag_of_distribution}
\mathcal{D}=\mathcal{D}^{(0)}\subset\mathcal{D}^{(1)}\subset\cdots\subset\mathcal{D}^{(\nu)}
\end{equation}the \emph{derived flag} of $D$.
Generically, this sequence ends up with the whole tangent bundle, that is, $\mathcal{D}^{(\nu)}=TM$. 
\begin{definition}\label{type:distribution}
Let $\mathcal{D}^{(\nu)}$ be as in \eqref{derived_flag_of_distribution}. If $\mathcal{D}^{(\nu)}$ has rank $\ell$, then $\mathcal{D}$ is said to be of \emph{type} $(\nu, \ell)$.
\end{definition}Notice that $\mathcal{D}^{(\nu)}$ is the smallest integrable distribution that contains $\mathcal{D}$ as a sub-bundle. Theorem~\ref{Frobenius} is the case that $\mathcal{D}$ has type $(0,k)$. 
Observe that if $\mathcal D$ has type $(\nu,\ell)$ there exists a non-degenerate system of $n-\ell $ first integrals since $\mathcal D^{(\nu)}$ is integrable in the sense of Frobenius.        This implies that  $M$ is foliated by  $(n-\ell)$-parameter family of invariant submanifolds of $\mathcal D$  (see Theorem \ref{THM2:INV_Controls}).    To discuss the existence of isolated (zero-parameter family) invariant submanifolds  we need the notion of the generalized first integral,  which was first introduced in \cite{AH}.

\begin{definition}\label{GFI}
Let $\mathcal{D}$ be a distribution. A non-degenerate set $\rho=(\rho^1,\ldots,\rho^d)$ of smooth real-valued functions is called a system of \emph{generalized first integrals} of $\mathcal{D}$ if \eqref{FI} holds on the zero locus of $\rho$, that is, for any section $X$ of $\mathcal{D}$,
\begin{equation}\label{generalized first}
(X\rho)(x)=0,\quad\text{for all}\;x\in M\quad\text{with}\;\rho(x)=0.
\end{equation}
\end{definition}In this paper, we focus our attention to the affine control system \eqref{affine_control_system} under an  assumption that the linear span of $f,g_1,\ldots,g_m$ has constant rank. 

Now we recall some basics of non-linear control systems. For other definitions and theorems we refer the readers to our basic references~\cite{AS, J97,NS,S73}. 
Consider   the set of points that are reachable from a point in some non-negative time.  
A point $q\in M$ is said to be \emph{reachable from} $p$ if 
$ q=\Gamma(p,t,u), $ for some $ u\in \mathcal U $  and for some  $t \ge 0.$   
 The \emph{reachable set $\mathcal{R}_{p}$ of the control system \eqref{affine_control_system}
 from a point $p\in M$} is a subset of $M$ defined by
\[
\mathcal{R}_{p}=\left\{\Gamma (p, t, u):  t\ge 0, ~~  u\in\mathcal U \right\}.
\]

The control system \eqref{affine_control_system} is said to be \emph{controllable from $p\in M$} if
\begin{equation}\label{defn:controllable_from_a_pt}
\mathcal{R}_{p}=M
\end{equation}and it is called \emph{controllable} if \eqref{defn:controllable_from_a_pt} holds for every $p\in M$.
We define the \emph{orbit $\mathcal{O}_{p}$ of the control system \eqref{affine_control_system} through a point $p\in M$} 
to be  the set of points $q\in M$
such that either $q$ is reachable from $p$ or $p$ is reachable from $q$.

The term $f(x)$ in (\ref{affine_control_system}) is called the \emph{drift}.  It is obvious that for  affine control systems  without drift   we have  
\begin{equation}\label{R-O relation}
\mathcal{R}_{p}=\mathcal{O}_{p}
\end{equation} for each  $p\in M.$
We observe also that the Nagano-Sussmann orbit theorem \cite{AS} holds for the affine control system  (\ref{affine_control_system}) regardless of the  drift, 
that is,  

\begin{theorem}\label{Sussman thm} $\mathcal{O}_p$ is a connected immersed submanifold of $M$.
\end{theorem}

Let $\mathcal G$    be   the Lie algebra of vector fields generated by $g_1, \ldots, g_m.$   
Consider  the vector space $\mathcal G(p) \subset T_pM, $   which  is the linear span at $p\in M$  of   the left iterated Lie brackets 
\[
[g_{i_{1}},[g_{i_{2}},\cdots,[g_{i_{k-1}},g_{i_{k}}]\cdots]]
\]of the  vector fields $g_1, \ldots, g_m.$   Then we have 
\begin{theorem}[Rashevsky-Chow theorem, \cite{AS}]\label{R-C:thm}
Let $M$ and $\mathcal{G}$ be as above. If $M$ is connected and ${\mathcal G }(p) =T_{p}M$ for a point $p\in M$, then $\mathcal{O}_{p}$ contains an open neighborhood of $p$. 
\end{theorem}
\vskip 1pc

Now  we state our main results as follows:  
\begin{theorem}\label{THM1:GFI}
Suppose that the vector fields $f,g_1,\ldots,g_m$ span a subspace $\mathcal{D}$ of constant rank. Let $N$ be a submanifold defined as the common zero set of real-valued functions $\rho=(\rho^1,\ldots,\rho^d)$ with $d\rho^1\wedge\cdots\wedge d\rho^d\neq0$. Then $N$ is invariant under controls of \eqref{affine_control_system} if and only if $\rho$ is a  system of generalized first integrals of $\mathcal{D}$.
\end{theorem}
\begin{theorem}\label{THM2:INV_Controls}
Let $\mathcal{D}$ be as in Theorem~\ref{THM1:GFI}. If $\mathcal{D}$ has type $(\nu, \ell)$,  $1\le \ell 
\le n, $  
then $M$ is foliated by $\ell$-dimensional  submanifolds that are invariant under controls of \eqref{affine_control_system}. This type condition is given by a system of partial differential equations of order $(\nu+1)$ for $f$ and $g_j$, $j=1,\ldots,m$. 
\end{theorem}

 Given a control system  proving  its controllability is primarily important in geometric control theory,  see \cite{I}.   In this paper we discuss its `non-controllability'  determining the existence of invariant submanifolds.   As a  possible application we would be able to  design  control systems with   prescribed invariant submanifolds.  
In \S\ref{method}  we construct the derived flag of $\mathcal D$ and find the generalized first integrals by using   various generalizations of the Frobenius theorem on involutivity  (cf.  \cite{I} and \cite{IKGM}).   \S\ref{exa} is devoted to examples of   invariant submanifolds.

\vskip 2pc
\section{Proof of the main results}\label{thms:pf}
\subsection{Proof of Theorem~\ref{THM1:GFI}}
Suppose that $\rho=(\rho^1, \cdots, \rho^d)$ is a non-degenerate system of real-valued functions whose common zero set $N$ 
 is invariant  under controls of \eqref{affine_control_system}.  To show that  $\rho$ is a  generalized first integral  it suffices to show that
\begin{equation}\label{assertion1}
\left(f\rho\right)(x)=\left(g_j\rho\right)(x)=0,\;j=1,\ldots ,m,\;\text{for all}\;x\in N, 
\end{equation}
since $f, g_1, \ldots, g_m$ span $\mathcal D.$
Let us fix a point $p\in N$. If we  choose  control $u(t)=0, $ to be constantly zero,  the trajectory $x(t)$
 with $x(0)=p$ is the integral curve of $f$ through $p$.   Since $N$ is invariant under controls of \eqref{affine_control_system},  we have by the chain rule 
 \begin{equation}\label{chain rule}
0=\frac{d}{dt}\rho(x(t))=\left(\dot{x}\rho\right)(x(t))=\left(f\rho\right)(x(t)).
\end{equation}
Substituting $t=0$ in \eqref{chain rule}, we obtain 
\begin{equation}\label{f_rho}
\left(f\rho\right)(p)=0.
\end{equation}  Now for each  $j=1,\ldots,m,$ let  $u(t)=(u_1(t),\ldots,u_{m}(t))\in\mathcal{U}$ defined by $u_{j}(t)=1$ and $u_{i}(t)=0$ for all $i\neq j$. Then there exists a unique trajectory $y(t)$ in some finite time interval with $y(0)=p$ which is the integral curve of 
\[
\dot{x} = f(x) + g_j(x).
\]
In a similar fashion to \eqref{chain rule} and \eqref{f_rho}, one can deduce that
\begin{equation}\label{f+g}
\left(\left(f+g_{j}\right)\rho\right)(p)=0.
\end{equation} 
From \eqref{f_rho} and \eqref{f+g} we have

\begin{equation}\label{gj_rho}
\left(g_j\rho\right)(p)=0.
\end{equation}
Since $p$ was arbitrary, \eqref{f_rho} and \eqref{gj_rho} prove the assertion \eqref{assertion1}.

Conversely, suppose that $\rho=(\rho^1,\ldots,\rho^d)$ is a non-degenerate system of generalized first integrals of $\mathcal{D}$. Let $N$ be the zero locus of $\rho$. We now fix a point $p\in N$.  Consider first  a constant control $u_c(t)= c, $ $c=(c_1, \ldots, c_m).$ Since $f$ and $g_j$ are tangent to $N$, we see that 
\[
f+\sum_{j=1}^m c_jg_j
\]
is tangent to $N$. This implies that  $\Gamma(p,t,u_c)\in N$, that is, $\rho(x(t))=0$. Now for a piecewise constant control $u=u(t),$ we apply the same argument to each piece to obtain  $\Gamma(p,t,u(t))\in N.$  Therefore, the orbit of $p$ is contained in $N$,  which implies that $N$ is invariant under the controls of \eqref{affine_control_system}.   \qed

\subsection{Proof of Theorem~\ref{THM2:INV_Controls}}
 Suppose that a smooth real-valued function $\rho$ is a first integral of $\mathcal{D}$ and  that $Z$ is 
 a  section of $[\mathcal{D},\mathcal{D}], $ namely,  $Z=[X,Y] $, for  some sections $X,Y$ of $\mathcal{D}$.
Then 
\[
Z\rho=(XY-YX)\rho=0.
\] This implies that $\rho$ is a first integral of  $\mathcal{D}^{(1)}$. By induction, $\rho$ is  a first integral of $\mathcal{D}^{(\nu)}$. Since $\mathcal{D}^{(\nu)}$ is integrable and has rank $\ell$, there exists a non-degenerate system $\rho=(\rho^1,\ldots,\rho^{n-\ell})$ of first integrals by Theorem~\ref{Frobenius}. Thus $M$ is foliated by integral manifolds $~\rho= constant~$ of $\mathcal{D}^{(\nu)}$, where the dimension of each leaf is $\ell$. Moreover, Theorem~\ref{THM1:GFI} implies that each leaf $N$ is an invariant submanifold under controls of \eqref{affine_control_system}.
Notice that the condition $[\mathcal{D}^{(\nu)},\mathcal{D}^{(\nu)}]\subset\mathcal{D}^{(\nu)}$ involves  the derivatives of $f,g_{j}$, $j=1,\ldots,m$, up to order $\nu+1$. This  will be shown more explicitly in the dual arguments of  \S\ref{method},  where  we construct  the derived flag of the associated Pfaffian system. \qed

\section{Construction of the derived flag and the generalized first integrals}\label{method}
This section is mainly concerned with the systematic construction of the derived flag and the generalized first integrals of a distribution, which are defined in \S\ref{intro}. For the sake of computation, we use the exterior differentiation and algebraic operations on differential ideals instead of vector fields and their brackets. This approach has  been used also  in \cite{AH,HP2,HP3} for the cases of real vector fields and in \cite{HK,HL,HP1} for the cases of complex vector fields. The notion of type for Pfaffian system (Definition~\ref{type:Pfaffian_system}) is due to R. B. Gardner~\cite{G}.

Let $M$ be a smooth real manifold of dimension $n$. Let $\Omega^{0}$ be the ring of smooth real-valued functions on $M$ and $\Omega^{k}$ $(1\leq k\leq n)$ the module over $\Omega^{0}$ of smooth $k$-forms on $M$ with smooth real coefficients. Then $\Omega^{*}:=\oplus_{k=0}^{n}\Omega^{k}$ is the exterior algebra equipped with operations of the wedge product $\wedge$ and the exterior differentiation $d$. Our standard reference for this section is \cite{BCGGG}. 
\begin{definition}
A subalgebra $\mathcal{I}$ of $\Omega^{*}$ is called an \emph{algebraic ideal} if 
\begin{itemize}
\item[(i)] $\mathcal{I}\wedge\Omega^{*}\subset\mathcal{I}$;
\item[(ii)] if $\phi=\sum_{k=0}^{n}\phi_{k}\in\mathcal{I}$, $\phi_{k}\in\Omega^{k}$, then each $\phi_{k}$ is in $\mathcal{I}$ (homogeneity condition).
\end{itemize}
\end{definition}
\noindent Note that the homogeneity condition implies that  $\mathcal{I}$ is a two-sided ideal, that is, 
\[
\Omega^*\wedge\mathcal{I}\subset\mathcal{I}.
\]Let
$\psi=(\psi^1,\ldots,\psi^l)$ be a system of smooth differential forms of any degree. We denote by $(\psi)$ the algebraic ideal generated by $\psi$. Then each element of $(\psi)$ has the form
\[
\sum_{k=1}^{l}\xi^k\wedge\psi^k
\]for some $\xi^k\in\Omega^{*}$. For $\alpha,\beta\in\Omega^{*}$,  if $\alpha-\beta\in(\psi)$  
we write  
$
\alpha\equiv\beta,\;\mathrm{mod}\;(\psi).$

Now let $X_1,\ldots, X_{p}$ be linearly independent smooth vector fields on $M$ and $\mathcal{D}$ a distribution generated by them. Consider a system of independent smooth 1-forms $\theta:=(\theta^1,\ldots,\theta^s)$, $s+p=n$, on $M$ which annihilates $X_1,\ldots,X_p$. We denote by $I$ the sub-bundle of the cotangent bundle $T^{*}M$ generated by $\theta$. Now we fix notations: For any sub-bundle $I$ of $T^{*}M$ we denote by $\underbar{I}$ the $\Omega^{0}$-module of smooth sections of $I$ and by $(I)$ the algebraic ideal of $\Omega^{*}$ generated by the smooth sections of $I$. Then the Frobenius integrability for $\mathcal{D}$ becomes
\begin{equation}\label{Frobenius_integrability:form_version}
d\underbar I\subset(I).
\end{equation}A set of real-valued functions $\rho=(\rho^1,\ldots,\rho^d)$ is a first integral if
\begin{equation}\label{FI:form_version}
d\rho\in\underbar I.
\end{equation}
\noindent Now we consider the composition of the exterior differentiation $d:\underbar{I}\rightarrow\Omega^{*}$ and the natural projection $\pi:\Omega^{*}\rightarrow\Omega^{*}/(I)$:
\[
\underbar I\xrightarrow{d}\Omega^{2}\xrightarrow{\pi}\Omega^{2}/(I).
\]
Then $\delta:=\pi\circ d$ is a module homomorphism. We define a sub-module $\underbar I^{(1)}$ of $\underbar I$ by setting
\[
\underbar I^{(1)}=\mathrm{ker}\;\delta.
\]      
Assuming that $\underbar I^{(1)}$ has constant rank, $I^{(1)}\subset T^{*}M$ is now a sub-bundle of ${T}^{*}M$ which we call the \emph{first derived system} of $I$. Then the sequence
\[
0\rightarrow \underbar I^{(1)}\rightarrow \underbar I\xrightarrow{\delta}d\underbar I/(I)\rightarrow0
\]is exact. Assuming that $\underbar I^{(k-1)}$ has constant rank on $M$, we define inductively the \emph{$k$-th derived system} $I^{(k)}$ of $I$ by the exactness of 
\[
0\rightarrow \underbar I^{(k)}\rightarrow \underbar I^{(k-1)}\xrightarrow{\delta} d\underbar I^{(k-1)}/(I^{(k-1)})\rightarrow0.
\]
 By setting  $\nu$  to be the smallest non-negative integer such that
\[
I^{(\nu)}=I^{(\nu+1)}  
\]we obtain a sequence
\[
I=I^{(0)}\supset I^{(1)}\supset\cdots\supset I^{(\nu-1)}\supset I^{(\nu)},
\]which we shall call the \emph{derived flag} of $I$. Notice that $I^{(\nu)}$ is the largest sub-bundle of $I$ that satisfies the Frobenius integrability \eqref{Frobenius_integrability:form_version}. Moreover, by using the formula 
\[
d\varphi(X,Y)=X\varphi(Y)-Y\varphi(X)-\varphi([X,Y])
\]for all $1$-form $\varphi$ and all smooth vector fields $X$ and $Y$, it is easy to see that for each $k=0,1,\ldots,\nu$,  
\begin{equation}\label{dual_relation}
{\mathcal{D}^{(k)}}^{\perp}=\underbar I^{(k)}.
\end{equation}

\begin{definition}\label{type:Pfaffian_system}
Let $I$ and $I^{(\nu)}$ be as above. If $I^{(\nu)}$ has rank $q$, then $I$ is said to be of \emph{type} $(\nu,q)$. 
\end{definition}Compare this to Definition~\ref{type:distribution} and observe that $I$ has type $(\nu,q)$ if and only if $\mathcal{D}$ has type $(\nu,n-q)$ by \eqref{dual_relation}. We  observe also that \eqref{FI:form_version} implies $d\rho\in I^{(1)}$ and inductively $d\rho\in I^{(\nu)}$. Thus we have the following proposition,  which is a basic observation of \cite{G}.

\begin{proposition}
A real-valued function $\rho$ is a first integral if and only if $d\rho\in I^{(\nu)}$. Therefore,  if $I$ has type $(\nu,q)$ then there is a $q$-parameter family of non-degenerate first integrals.
\end{proposition}
Now the defining property \eqref{generalized first} of the generalized first integral 
states as follows:
\begin{proposition}  On a smooth manifold $M^n$,  let  $\mathcal{D}$  be a distribution spanned by vector fields $X_1, \ldots, X_p, $ and $\theta =(\theta^1, \ldots, \theta^s), $
$p+s=n, $ be $1$-forms that annihilate $\mathcal D.$   Let  $\rho=(\rho^1,\ldots,\rho^d)$, $d\leq s$, be  a non-degenerate set of smooth real-valued functions.  Then $\rho$ 
 is  a system of generalized first integrals of $\mathcal{D}$ if and only if 
\begin{equation}\label{GFI:Form}
d\rho^\mu\in(\rho, \theta),\;\mu=1,\ldots,d.
\end{equation} 
\end{proposition}
Now let $\omega^1,\ldots,\omega^p$ be a set of closed smooth $1$-forms that completes $\theta$ to a local coframe
\begin{equation}\label{coframe}
(\theta^1,\ldots,\theta^s,\omega^1,\ldots,\omega^p)
\end{equation}of $M$. Set
\begin{equation}\label{Torsion}
d\theta^l\equiv\sum_{j<k}T^{l}_{j,k}\omega^j\wedge\omega^k,\;\mathrm{mod}\;(\theta),\;l=1,\ldots,s.
\end{equation}In matrices, \eqref{Torsion} is  written as 
\begin{equation}\label{torsion:original}
\begin{bmatrix}
d\theta^1\\
\vdots\\
d\theta^s
\end{bmatrix}
\equiv
\underbrace{
\begin{bmatrix}
T^{1}_{1,2}& T^{1}_{1,3}&\cdots& T^{1}_{p-1,p}\\
\vdots&&&\vdots\\
T^{s}_{1,2}& T^{s}_{1,3}&\cdots& T^{s}_{p-1,p}
\end{bmatrix}
}_{\mathcal{T}}
\begin{bmatrix}
\omega^1\wedge\omega^2\\
\omega^1\wedge\omega^3\\
\vdots\\
\omega^{p-1}\wedge\omega^{p}
\end{bmatrix},\;\mathrm{mod}\;(\theta).
\end{equation}    $\mathcal{T}$ in \eqref{torsion:original} is a matrix of size $s\times  {p\choose 2}$, called the \emph{torsion matrix} with respect to the coframe \eqref{coframe}. If $\mathcal{T}\equiv0$, this is the case of the Frobenius theorem. If $\mathcal{T}\notequiv0$, then we find generators of $I^{(\nu)}$ as follows: Suppose that 
\begin{equation}\label{how_to_find_generators_of_the_first_derived_System}
\phi:=\sum_{\gamma=1}^{s}a_{\gamma}\theta^{\gamma}\in I^{(1)}=\mathrm{ker}\;\delta.
\end{equation}Then we have 
\begin{equation}\label{torsion_rank_assertion1}
d\phi\equiv\sum_{\gamma=1}^{s}a_{\gamma}d\theta^{\gamma}\equiv0,\;\mathrm{mod}\;(\theta).
\end{equation}Substituting \eqref{Torsion}  for  $d\theta^{\gamma}$ in \eqref{torsion_rank_assertion1}, it follows from the independence of $\omega^k\wedge\omega^l$ in $\Omega^{2}/(I)$ that the row vector $(a_1,\ldots,a_s)$ belongs to the left null-space of $\mathcal{T}$. By  finding the generators of the left null-space of $\mathcal{T}$ we obtain a set of generators \eqref{how_to_find_generators_of_the_first_derived_System} of $I^{(1)}$. Then by the linear algebra of  the torsion matrix for the Pfaffian system $I^{(1)}$, we obtain the generators of $I^{(2)}$, and for $I^{(\nu)}$ by induction.

We use similar arguments for the generalized first integrals. Suppose that $\rho=(\rho^1,\ldots,\rho^d)$, $d\leq s$, is a non-degenerate system of generalized first integrals of $\mathcal{D}$. Then, for each $\mu$, $1\leq \mu\leq d$, $d\rho^\mu\in(\rho,\theta)$, namely
\begin{equation}\label{GFI:assertion1}
d\rho^\mu=\sum_{\lambda=1}^{d}\rho^\lambda\psi^{\mu}_{\lambda}+\sum_{\gamma=1}^{s}b^{\mu}_{\gamma}\theta^{\gamma}
\end{equation}for some smooth $1$-forms $\psi^\mu_\lambda$ and some smooth functions $b^\mu_\gamma$. Applying $d$ to \eqref{GFI:assertion1}, we have 
\begin{equation}\label{GFI:assertion2}
0\equiv\sum_{\gamma=1}^{s}b^{\mu}_{\gamma}d\theta^{\gamma},\;\mathrm{mod}\;(\rho,\theta),\;\mu=1,\ldots,d.
\end{equation}Consider the matrix of coefficients
\[
\mathcal{B}=\begin{bmatrix}
b^{1}_{1}&\cdots& b^{1}_{s}\\
\vdots& &\vdots\\
b^{d}_{1}&\cdots &b^{d}_{s}
\end{bmatrix}.
\]  
In matrices \eqref{GFI:assertion2} can be written as 
\[
\begin{bmatrix}
0\\
0\\
\vdots\\
0
\end{bmatrix}
\equiv
\mathcal{B}\mathcal{T}
\begin{bmatrix}
\omega^1\wedge\omega^2\\
\omega^1\wedge\omega^3\\
\vdots\\
\omega^{p-1}\wedge\omega^{p}
\end{bmatrix},\;\mathrm{mod}\;(\rho,\theta).
\]
\noindent Since  $\omega^k\wedge\omega^l$ are independent $2$-forms we have
\[
\mathcal{B}\mathcal{T}\equiv0,\;\mathrm{mod}\;(\rho).
\]Moreover,  since $(\rho^1, \ldots, \rho^d)$ are non-degenerate   $\mathcal{B}$ has maximal rank $d$ on the zero locus of $\rho,  $  which implies that   $\mathrm{rank}\;\mathcal{T}\leq s-d$ on the zero locus of $\rho$.  
Thus  we observed that \textit{any square sub-matrix of $\mathcal{T}$ of size $s-d+1$ has determinant zero on the zero locus of $\rho$.} Then  we factor those determinants,  find non-degenerate functions that generate all those determinants, which are the candidates to be the generalized first integral $\rho.$ Finally, if thus found $\rho$ satisfies  \eqref{GFI:Form}, then $\rho$ is the desired system of generalized first integrals. 

\section{Examples}\label{exa}  

The examples we present in this section are rather artificial,  devised simply to show how our method works.
\begin{example}[a single invariant submanifold] Consider an affine control system on $\mathbb R^3=\{(x,y,z)\}$ given by 
\begin{equation}\label{example 1}
\begin{bmatrix}
\dot{x}\\
\dot{y}\\
\dot{z}
\end{bmatrix}
=
\underbrace{\begin{bmatrix}
1\\
y\\
0
\end{bmatrix}}_{g_1}
u_{1}
+
\underbrace{\begin{bmatrix}
0\\
1\\
xz
\end{bmatrix}}_{g_2} u_{2}:=g_1u_1+g_2u_2,
\end{equation}
where $\mathbb{R}^3=\{(x,y,z)\}$ is the state space and  $(u_1,u_2)\in\mathbb{R}^2$ are controls. Then
\[\theta=xyz ~dx-xz~dy+dz
\]annihilates  
$g_1$ and $g_2$. We take up a coframe
$\{\theta,dx,dy\}$. Then with respect to this coframe 
\[
d\theta\equiv
\underbrace{
-z(1+x)}_{\mathcal T}
  dx\wedge dy,\;\mathrm{mod}\;(\theta).
\]
In this case the torsion is a $1\times 1$ matrix $\mathcal{T}=-z(1+x).$
Any non-degenerate factor of $T$, in particular, $\rho(x,y,z)=z$ is a candidate for generalized first integral. Now we shall examine the condition \eqref{GFI:Form} for $\rho(x,y,z)=z$. Indeed, we have
\[
d\rho=dz=\theta-\rho(xy dx-xdy)\in(\rho, \theta).
\]     
Therefore, $\{z=0\}$ is an invariant submanifold.  Since $g_1$ and $g_2$ are independent on 
$z=0$ the system \eqref{example 1} restricted on $z=0$ is controllable.
\end{example}

\medskip

\begin{example}[no invariant submanifolds]
We slightly change \eqref{example 1} to 
\[
\begin{bmatrix}
\dot{x}\\
\dot{y}\\
\dot{z}
\end{bmatrix}
=
\underbrace{\begin{bmatrix}
1\\
y\\
0
\end{bmatrix}}_{h_1}
u_{1}
+
\underbrace{\begin{bmatrix}
0\\
1\\
xy
\end{bmatrix}}_{h_2} u_{2}:=h_1u_1+h_2u_2.
\]
 In the same way as in the previous example we have
\[
\theta=xy^2dx - xydy + dz,
\]so that
\[
d\theta=\underbrace{   -y(1+2x)}_{\text{torsion}} dx\wedge dy.
\]
Now $\rho(x,y,z)=y$ is a non-degenerate factor of the torsion,  which is unique modulo multiplication by unit. Hence the only candidate to be an invariant submanifold is $y=0$.   However, 
\[
d\rho = dy \notin (y, \theta),
\]therefore, there are no invariant submanifolds.
\end{example}
 
\medskip
 
\begin{example}[foliation by invariant submanifolds, vehicle on a slanted plane]  First we describe the motion of a car on a slanted plane as follows: 
Let $M=\{(x,y,z,w): x,y,z\in\mathbb{R}^3,   w\in S^1 \} $ be the state space of the  affine control system
\[
\begin{bmatrix}
\dot{x}\\
\dot{y}\\
\dot{z}\\
\dot{w}\\
\end{bmatrix}
=
\begin{bmatrix}
\cos\phi_0\cos w\\
\sin w\\
\sin \phi_0 \cos w\\
0
\end{bmatrix}u_{1}+
\begin{bmatrix}
0\\
0\\
0\\
1
\end{bmatrix}u_{2},
\]where $(x,y,z)$ is the position of the center of mass of the car and $w$ is 
the angle measured from the first coordinate axis of the slanted plane to the direction of the vehicle.  Assuming the slanted plane meets with the $xy$-plane along $y$-axis,  the first coordinate line of the slanted plane  intersects with $y$-axis perpendicularly,   and 
 $\phi_{0}, $ $0<\phi_0<\pi/2,$  is the angle  between  these two planes.  We assume  two possible motions of a car: one can drive the car forward and backwards with a fixed linear velocity $u_1\in\mathbb{R}$, and turn the car around its center of mass with a fixed angular velocity $u_{2}\in\mathbb{R}$. We note that this control system is a modification of the typical model of a car in $\mathbb{R}^2\times S^1$~(cf.~\cite[Example $1.29$]{AS}). However, we shall work in a more general setting so that $a$ and $b$ are assumed to be any positive real numbers. We consider
 
\begin{equation}\label{car_model:slanted}
\underbrace{\begin{bmatrix}
\dot{x}\\
\dot{y}\\
\dot{z}\\
\dot{w}\\
\end{bmatrix}}_{\dot{\mathbf x}}
=
\underbrace{\begin{bmatrix}
a\cos w\\
\sin w\\
b \cos w\\
0
\end{bmatrix}}_{g_1} u_{1}+
\underbrace{\begin{bmatrix}
0\\
0\\
0\\
1
\end{bmatrix}}_{g_2}u_{2}:=g_{1}u_{1}+g_{2}u_{2}.
\end{equation}
\noindent The system of $1$-forms $\theta=(\theta^1, \theta^2)$ given by 

\begin{equation*}
\begin{aligned}
\theta^1&= b ~dx - a~ dz,\\
\theta^2&=b \cos w ~dy -  \sin w~ dz
\end{aligned}
\end{equation*}   
\noindent annihilates  $g_1$ and $g_2$.  Let    $I:= (\theta)$ be the ideal.
Then   we have
\begin{equation*}
\begin{aligned}\begin{bmatrix}
d\theta^1\\
d\theta^2
\end{bmatrix}
&=
\begin{bmatrix}
0\\
-b \sin w ~dw\wedge dy -  \cos w ~dw\wedge dz
\end{bmatrix}\\ &
\equiv \begin{bmatrix}
0\\
 \sec w~ dz\wedge dw\end{bmatrix},  ~~~ \;\mathrm{mod}\;(\theta^2). \end{aligned}
\end{equation*}

\noindent We see that  $I$ has type $(1,1)$  and the first derived system $I^{(1)}$  is generated by $\theta^1.$ 
 The first integral is 

\[
\rho(x,y,z,w)= bx - az,
\] and therefore, 
hyperplanes 
\begin{equation}\label{slant inv}
bx-az=\text{constant}
\end{equation}are invariant under controls of \eqref{car_model:slanted}. Since  $g_1, g_2 $ and 
\[
[g_1, g_2]=(a\sin w, -\cos w, b\sin w, 0)^t
\]are independent on  each hyperplane of \eqref{slant inv} the control system \eqref{car_model:slanted} restricted to each of these invariant hyperplanes is controllable by \eqref{R-O relation} and Theorem~\ref{R-C:thm}.
\end{example}

\medskip

\begin{example}[slanted plane with drift]
Consider an affine control 
\begin{equation}\label{slant with drift}
\dot{\mathbf x} = f({\mathbf x})+ g_1({\mathbf x}) u_1 + g_2({\mathbf x}) u_2, \end{equation}
where $\mathbf x$ and $g_1, g_2$ are the same as in \eqref{car_model:slanted} and $f$ is a drift.
The following are obvious:
\vskip 1pc
\noindent a) If $f$ is contained in the linear span of ($g_1, g_2)$ then $\rho(x,y,z,w)=bx-az$ is a first integral,  and therefore,  $M$ is foliated by the invariant submanifolds $\rho=\mathrm{constant}.$

\noindent b) If $f$ has a non-zero component only in $y$-direction,  then $\rho(x,y,z,w)=bx-az$ is a first integral because $\rho$ is independent of $y$ variable.

\noindent c) For  generic choices of $f$    invariant submanifolds do not exist.
\vskip 1pc

A single invariant submanifold, or equivalently, a generalized first integral, occurs only under 
special  assumptions on $f$.    For instance,   if $f = (\phi, 1, 0, 0)^t, $  where $\phi$ is any function that is divisible by $\rho=bx-az, $   then $ \rho $ is a generalized first integral, therefore, a single hyperplane $\rho=0$ is invariant under controls of \eqref{slant with drift}. This is obvious because $f=(0,1,0,0)^t$ on $\rho=0$, so that $f$ has only $y$-directional component.   In fact, 
\begin{equation}\label{example4-1}
\theta = -b \cos w ~dx + b \phi \cos w ~ dy + (a\cos w - \phi\sin w) dz 
\end{equation}
annihilates $f, g_1$ and $g_2.$ 
Since $d\rho=b dx - a dz $   and $\phi$ is divisible by $\rho,$  rearranging the right hand side of \eqref{example4-1} we have 
\begin{equation}\label{last example1} \theta = -\cos w ~d\rho + \rho\alpha (b\cos w ~dy - \sin w ~dz), \quad\text{for some function }\alpha.
\end{equation}   Solving  \eqref{last example1} for $d\rho$ yields  
 \eqref{GFI:Form}.
Observe also that  
\begin{equation}\label{last example2} d\theta \equiv -\phi \sec w~ dw\wedge dz + b\cos w~d\phi\wedge dy - \sin w ~d\phi\wedge dz, \quad \text{ mod } (\theta). \end{equation}
Since $\phi=\rho\alpha, $ for some  $\alpha$, \eqref{last example2} implies that 
\begin{equation}\label{last 2} d\theta = \rho\wedge\Psi^1 + d\rho\wedge \Psi^2 + \theta\wedge\Psi^3, \quad\text{for some } \Psi^j\in \Omega^*, j=1,2,3.\end{equation}
For the tangent vectors to $\rho=0$ the first two terms of the right hand side are zero and \eqref{last 2} implies that $\theta$ is Frobenius integrable on $\rho=0.$   But in this particular case, the foliation  has a single leaf: the submanifold $\rho=0$ itself. 
\end{example}

\subsection*{Acknowledgements}
The second author was supported by National Research Foundation of Korea with grant NRF-2015R1A2A2A11001367.


\end{document}